# STRUCTURAL PROPERTIES OF PROPORTIONAL FAIRNESS: STABILITY AND INSENSITIVITY

By Laurent Massoulié

*Thomson Paris Research Lab*

In this article we provide a novel characterization of the proportionally fair bandwidth allocation of network capacities, in terms of the Fenchel–Legendre transform of the network capacity region. We use this characterization to prove stability (i.e., ergodicity) of network dynamics under proportionally fair sharing, by exhibiting a suitable Lyapunov function. Our stability result extends previously known results to a more general model including Markovian users routing. In particular, it implies that the stability condition previously known under exponential service time distributions remains valid under so-called phase-type service time distributions.

We then exhibit a modification of proportional fairness, which coincides with it in some asymptotic sense, is reversible (and thus insensitive), and has explicit stationary distribution. Finally we show that the stationary distributions under modified proportional fairness and balanced fairness, a sharing criterion proposed because of its insensitivity properties, admit the same large deviations characteristics.

These results show that proportional fairness is an attractive bandwidth allocation criterion, combining the desirable properties of ease of implementation with performance and insensitivity.

**1. Introduction.** The abstract network bandwidth allocation (NBA) problem can be formulated as follows. A network supports connections of distinct types, indexed by $r$, the index $r$ spanning the set of types $\mathcal{R}$, assumed finite. Given the number $x_r$ of users of each type $r \in \mathcal{R}$, with $x_r \in \mathbb{N}$, the problem is to determine the total capacity allocated to type $r$ users, denoted be $\lambda_r$, with $\lambda_r \in \mathbb{R}_+$. The quantity $\lambda_r$ represents the rate at which data is received collectively by all users of type $r$. The allocation vector $\lambda := \{\lambda_r\}_{r \in \mathcal{R}}$ is constrained to lie in a set $\mathcal{C} \subset \mathbb{R}_+^{|\mathcal{R}|}$.

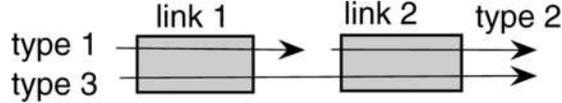

Fig. 1. *Example of a two-link network supporting three types of users.*

The set $\mathcal{C}$ is a suitable abstraction of all the physical capacity constraints of the actual network under consideration. An example of a two-link network is represented in Figure 1. This network supports three types of users, data destined to type-1 users going through link 1 only, data to type-2 users going through link 2 only, while data for type-3 users goes through the two links. Thus, when the two links have unit capacities, the corresponding set $\mathcal{C}$ is given by $\{\lambda \in \mathbb{R}_+^3 : \lambda_i + \lambda_3 \leq 1,\ i = 1, 2\}$. This example can be extended to the case where the network consists of an arbitrary number of links, and user types $r$ are characterized by collections of links used by data destined to them. Denoting by $\mathcal{L}$ the collection of links, and by $c_\ell$ the capacity of link $\ell \in \mathcal{L}$, the corresponding network capacity region then takes the form

$$\mathcal{C} = \left\{\lambda \in \mathbb{R}_+^\mathcal{R} : \sum_{r \in \mathcal{R}} A_{\ell r} \lambda_r \leq c_\ell, \ell \in \mathcal{L}\right\},$$

where $A_{\ell r}$ equals 1 or 0 according to whether type-$r$ users require capacity at link $\ell$ or not. Such types of capacity constraints have been considered for instance in [13, 14], as suitable models of wired networks with fixed routing such as the Internet, the matrix $A$ then reflecting the route that data of users of given type follows through the network. More general polyhedral capacity sets $\mathcal{C}$ arise when users of a given type $r$ can send data along several distinct routes through the network. Yet more general, nonpolyhedral (albeit still convex) capacity sets can adequately model the impact of interferences between data transmissions of distinct types in wireless networks; see [4] for such examples.

In the present work we only require the set $\mathcal{C}$ to be convex, and nonincreasing, that is to say, for any two vectors $\lambda, \lambda' \in \mathbb{R}_+^\mathcal{R}$ such that $\lambda_r \leq \lambda'_r, r \in \mathcal{R}$, then $\lambda$ belongs to $\mathcal{C}$ whenever $\lambda'$ does. These two assumptions are met in all the examples above mentioned.

Mo and Walrand [19] introduced the following criterion for determining the allocation vector $\lambda$. Given weights $w_r > 0$, and a parameter $\alpha \geq 0$, the so-called $(w, \alpha)$-fair allocation vector is the solution to the optimization problem

$$\max_{\lambda \in \mathcal{C}} \sum_{r \in \mathcal{R}} x_r U_r(\lambda_r / x_r), \tag{1}$$

where

(2) $$U_r(y) = \begin{cases} \dfrac{w_r}{1-\alpha} y^{1-\alpha}, & \text{if } \alpha \neq 1, \\ w_r \log(y), & \text{if } \alpha = 1. \end{cases}$$

This parametric family of allocation criteria contains the so-called proportional fairness criterion, introduced by Kelly [13], which corresponds to the special case $\alpha = 1$ and $w_r \equiv 1$. In the limit $\alpha \to \infty$, the $(w, \alpha)$-fair allocation coincides with the so-called max-min fair allocation (see [2] for a definition).

The rationale for proportional fairness, as explained in [13], lies in the following desirable decomposition property. Assume that the ultimate goal of bandwidth allocation is to maximize the sum of utility functions, $U_r$, of the rates $\lambda_r/x_r$ allocated to users of class $r$, exactly as in equation (1), these utility functions being known to the users but not to the network. Then the decomposition result of [13] states that this can be done by letting on the one hand the network allocate bandwidth according to proportional fairness, with weights $w_r$ specified by the network users, and on the other hand the network users selecting these weights $w_r$ appropriately, given their (privately known) utility functions $U_r$, and the network allocation in response to distinct weights $w_r$.

Alternatively, the unweighted proportional fairness allocation arises naturally from results in bargaining theory, in contrast to the above justifications based on microeconomic theory. Indeed, the results of Stefanescu and Stefanescu [27] (see also [17] for further discussion) imply that it is the only allocation of bandwidth that satisfies four natural axioms introduced by Nash [20] (namely, invariance with respect to affine utility transformations, Pareto optimality, independence of irrelevant alternatives, and symmetry), assuming that users' utility is a linear function of the rate they receive. (Note the difference with the previous microeconomic framework, which allowed arbitrary concave utility functions.) It is in fact the natural extension of Nash's bargaining solution, originally derived in the special context of two users, to an arbitrary number of users.

The rationales for candidate NBA solutions we have just reviewed originate from microeconomic theory of utility, and game (bargaining) theory, and assume a static set of network users. There is another line of approach to the NBA problem, which is essentially motivated by performance issues in a dynamic setting.

Specifically, assume that network users arrive and leave the system, the arrivals of type $r$ users being at the instants of a Poisson process of rate $\nu_r$. Assume further that users remain in the system until they have transferred a file of a given size, files associated with type $r$ users being exponentially distributed with parameter $\mu_r$. The state variable $x = \{x_r\}_{r \in \mathcal{R}}$ is then a

Markov process, with nonzero transition rates

(3)
$$x_r \to x_r + 1 \quad \text{with rate } \nu_r,$$
$$x_r \to x_r - 1 \quad \text{with rate } \mu_r \lambda_r.$$

A suitable rationale for selecting a NBA is to guarantee desirable properties of the above Markov process. One such property is stability (or equivalently, ergodicity), as it in turn implies that sojourn times of users are almost surely finite. Ergodicity cannot be guaranteed for all sets of traffic parameters $\nu_r$, $\mu_r$ and network capacity sets $\mathcal{C}$. In particular, letting $\rho_r := \nu_r/\mu_r$ denote the load brought by type $r$-users, when the vector $\rho = \{\rho_r\}_{r \in \mathcal{R}}$ does not belong to the capacity set $\mathcal{C}$, the process cannot be ergodic (for a proof, see, e.g., [4]). When $\rho$ is on the boundary of $\mathcal{C}$, Kelly and Williams [15] have established that the process cannot be positive recurrent, for sets $\mathcal{C}$ corresponding to wired networks with fixed routes. Their proof extends to the case of general convex nonincreasing capacity sets $\mathcal{C}$ with minor modifications.

A reasonable performance requirement is thus that, provided the traffic intensity vector $\rho$ lies in the interior $\overset{\circ}{\mathcal{C}}$ of $\mathcal{C}$, then the above Markov process is ergodic. Such a property is in fact satisfied for all $(w, \alpha)$-fair bandwidth allocation criteria, as follows from the Lyapunov function-based stability proof of Bonald and Massoulié [3] (see also de Veciana, Lee and Konstantopoulos [10] who first established the result for the case of max-min fairness, Ye [29] and Key and Massoulié [16] for an extension to more general utility functions $U_r$ in the allocation definition (1)).

Thus, the requirement of achieving ergodicity for the largest possible set of traffic intensity vectors $\rho$, being met by all $(w, \alpha)$-fair NBA, does not distinguish one such criterion as superior to the others.

A more stringent requirement has been suggested by Bonald and Proutière [5], namely that not only the stability region (defined to be the set of vectors $\rho$ such that the system is ergodic) be maximal, but also that the corresponding Markov process be *insensitive* to the distribution of sizes of the files transferred by each class of users. Roughly speaking, insensitivity means that the stationary distribution of the numbers of users in the system is unaffected if the service time distributions are modified, provided their mean is left unchanged. For characterizations of insensitive systems, we refer the reader to Schassberger [25] and references therein. In particular, it holds that, when service rate to users of one type is shared equally among such users, that is to say, under a processor sharing assumption, reversibility of the original Markov process ensures it is insensitive [5, 25].

If insensitivity holds, the system remains ergodic under the natural stability condition for arbitrary phase-type (i.e., mixtures of convolutions of exponential distributions; see, [1], page 80), not necessarily exponential, service time distributions. Note that ergodicity under the natural stability conditions, and for general, nonexponential service time distributions, has so far been established for the max-min fairness NBA in a recent article of Bramson [6], but a similar result has been missing for all other $(w, \alpha)$ NBAs. Although restricted to max-min fairness, the results of Bramson apply under very weak integrability assumptions on the service time distributions, and are not restricted to phase-type distributions.

Bonald and Proutière have identified a new NBA, the so-called *balanced fairness allocation*, which meets the two requirements of maximal stability region and insensitivity, and moreover maximizes the fraction of time during which the system is empty, among all allocations meeting these two requirements. They have furthermore identified special network topologies for which balanced fairness coincides with proportional fairness, and have shown that for all other network topologies, balanced fairness is distinct from any utility maximization NBA.

This leaves several questions open regarding the choice of an NBA. On the one hand, utility maximization allocations, such as $(w, \alpha)$-fairness or more specifically proportional fairness, can be implemented in a distributed manner (see, e.g., the seminal paper by Kelly, Maulloo and Tan [14]), and are motivated by microeconomic theory and game theory arguments in a static setting. In addition, they satisfy the criterion of maximal stability region in the dynamic setting, but do not seem to meet the more stringent requirement of insensitivity. On the other hand, balanced fairness does meet the latter requirement, but no simple distributed technique for realizing this NBA is known, if we except the special network topologies, identified in [5], where it coincides with proportional fairness.

In the present work, we provide a novel characterization of proportional fairness, and use it to improve upon this unsatisfactory state of affairs. Indeed, relying on this structural property, we show that the seemingly fortuitous coincidence of balanced fairness and proportional fairness on specific network topologies in fact reflects a deeper relationship between the two NBAs, that holds for any network topology as captured by the set $\mathcal{C}$. More precisely, we exhibit a third NBA, namely modified proportional fairness, which coincides in some asymptotic sense with proportional fairness. Under modified proportional fairness, the system is reversible, and hence insensitive. Furthermore, the steady state distributions under modified proportional fairness and balanced fairness admit the same large deviations characteristics, described by a simple explicit rate function.

As a by-product, we give a new proof of ergodicity of proportional fairness, which extends to a more general model of network dynamics including Markovian users routing. This in turn implies that the usual stability conditions still hold with service time distributions that are of phase type rather than exponential.

In view of these results, proportional fairness is an attractive candidate as a default NBA. Indeed, it is motivated by the following factors: (1) the decomposition property of [13], (2) axiomatic arguments from bargaining theory [27], (3) as an implementable approximation to balanced fairness, meeting the additional criteria of performance and insensitivity.

The structure of the paper is as follows. Section 2 gives the novel characterization of proportional fairness. Stability properties with Markovian user routing are proven in Section 3. The special case of phase type service distributions is discussed in Section 4. Section 5 establishes the relationships between balanced fairness and modified proportional fairness, and in particular the fact that the corresponding equilibrium distributions have the same large deviations characteristics.

**2. Characterization of proportional fairness via convex duality.** It is convenient to consider the logarithms of the allocated capacities $\lambda_r$, rather than the $\lambda_r$ themselves. Denote by $K$ the subset of $\mathbb{R}^{|\mathcal{R}|}$ in which these must lie, that is,

$$\gamma = \{\gamma_r\} \in K \Leftrightarrow \lambda = \{\exp(\gamma_r)\} \in \mathcal{C}.$$

Given $\gamma$, $\gamma'$ in $K$, and $\varepsilon \in (0,1)$, by convexity of the exponential function, for all $r \in \mathcal{R}$, one has

$$\exp(\varepsilon \gamma_r + (1-\varepsilon)\gamma'_r) \leq \varepsilon \exp(\gamma_r) + (1-\varepsilon)\exp(\gamma'_r),$$

and thus since $\mathcal{C}$ is convex, nonincreasing, then so is $K$. Denote by $\gamma^{\mathrm{PF}}(x)$ the vector of logarithms of proportionally fair allocations, that necessarily belong to $K$.

Denote by $\delta_K$ the function that equals zero on $K$, and $+\infty$ outside of $K$. The original characterization of $\lambda^{\mathrm{PF}}(x)$ as a maximizer of $\sum_{r \in \mathcal{R}} x_r \log(\lambda_r)$ over $\lambda \in \mathcal{C}$ readily implies that

(4) $$\gamma^{\mathrm{PF}}(x) \in \arg\sup_{\gamma \in \mathbb{R}^{\mathcal{R}}} (\langle \gamma, x \rangle - \delta_K(\gamma)).$$

Let now $\delta_K^*$ denote the Fenchel–Legendre convex conjugate function of $\delta_K$, that is,

$$\delta_K^*(x) = \sup_{\gamma \in \mathbb{R}^{\mathcal{R}}} (\langle \gamma, x \rangle - \delta_K(\gamma)).$$

Recall that the subgradient of a convex function $J$ defined on $\mathbb{R}^n$ at a point $x \in \mathbb{R}^n$, which is denoted by $\partial J(x)$, is the set consisting of all the vectors $h$ such that, for all $y \in \mathbb{R}^n$,

$$J(x) + \langle h, y - x \rangle \leq J(y).$$

We then have the following compact characterization of the function $\gamma^{\mathrm{PF}}$.

LEMMA 1. *The function $\gamma^{\mathrm{PF}}$ satisfies for all $x \in \mathbb{R}_+^{\mathcal{R}}$,*

$$\gamma^{\mathrm{PF}}(x) \in \partial \delta_K^*(x), \tag{5}$$

*where $\partial \delta_K^*(x)$ denotes the subgradient of the convex function $\delta_K^*$ at $x$.*

PROOF. It follows from Theorem 23.5, page 218 in Rockafellar [22] that conditions (4) and (5) are equivalent for any *proper* convex function $\delta_K^*$. Recall that a convex function is proper if it nowhere takes the value $-\infty$, and it takes finite values at some points. Both conditions hold for $\delta_K^*$, which establishes the lemma. □

This simple result allows to use the powerful theory of convex duality in the study of the function $x \to \gamma^{\mathrm{PF}}(x)$. For instance, we have the following:

LEMMA 2. *The function $\delta_K^*$ is continuously differentiable on $(0, \infty)^{\mathcal{R}}$, and thus on $(0, \infty)^{\mathcal{R}}$, $\gamma^{\mathrm{PF}}(x)$ coincides with the ordinary gradient of $\delta_K^*$ at $x$, and depends continuously on $x$.*

PROOF. By Theorem 25.1, page 242 in [22], at a point $x$ where a convex function admits a unique subgradient, it is differentiable, and its subgradient reduces to its ordinary gradient. The original allocation vector $\lambda^{\mathrm{PF}}(x)$ is uniquely defined at $x$ whenever $x_r > 0$ for all $r \in \mathcal{R}$, by strict concavity of the log function. Thus, $\gamma^{\mathrm{PF}}(x)$ is also uniquely defined at $x \in (0, \infty)^{\mathcal{R}}$, and hence it coincides with the ordinary gradient of $\delta_K^*$ at $x$.

Furthermore, by Theorem 25.5, page 246 in [22], the gradient of a proper convex function is continuous on the domain where the function is differentiable. The claimed continuity of the allocation vector $\gamma^{\mathrm{PF}}(x)$ on $x \in (0, \infty)^{\mathcal{R}}$ follows. □

Introduce now the alternative NBA, denoted PF′ for modified proportional fairness, and defined by

$$\lambda_r^{\mathrm{PF}'}(x) = \begin{cases} \exp(\delta_K^*(x) - \delta_K^*(x - e_r)), & \text{if } x_r > 0, \\ 0, & \text{otherwise.} \end{cases}$$

Define the function $L$ on $\mathbb{R}_+^{\mathcal{R}}$ by

$$L(x) = \delta_K^*(x) - \sum_{r \in \mathcal{R}} \log(\rho_r) x_r, \tag{6}$$

where $\rho_r = \nu_r/\mu_r$, $r \in \mathcal{R}$. It is readily verified that, under the PF′ allocation strategy, the Markov process is reversible, and thus insensitive. Indeed, one easily shows that the measure $\pi^{\mathrm{PF}'}$ on $\mathbb{Z}_+^{\mathcal{R}}$ by

$$\pi^{\mathrm{PF}'}(x) = \exp(-L(x)) \tag{7}$$

verifies the detailed balance equations

$$\pi^{\mathrm{PF}'}(x+e_r)\mu_r\lambda_r^{\mathrm{PF}'}(x+e_r) = \pi^{\mathrm{PF}'}(x)\nu_r, \qquad r \in \mathcal{R}, x \in \mathbb{Z}_+^{\mathcal{R}}.$$

The natural stability condition is, as discussed previously, the following:

$$\rho \in \overset{\circ}{\mathcal{C}}. \tag{8}$$

The following lemma gives useful properties satisfied by function $L$:

LEMMA 3. *The function $L$ is lower semicontinuous on $\mathbb{R}^{\mathcal{R}}$, and continuous on $\mathbb{R}_+^{\mathcal{R}}$. Furthermore, under assumption* (8), *there exist positive constants $a, A > 0$ such that for all $x \in \mathbb{R}_+^{\mathcal{R}}$,*

$$a\|x\|_\infty \leq L(x) \leq A\|x\|_\infty, \tag{9}$$

*where $\|x\|_\infty := \sup_{r \in \mathcal{R}} |x_r|$.*

PROOF. The function $\delta_K^*$ is lower semicontinuous, as the Fenchel–Legendre conjugate of a proper convex function (by Theorem 12.2, page 104 in [22]). The sum of an affine—and hence continuous—function with a lower semicontinuous function is lower semicontinuous. Thus $L$ is lower semicontinuous.

Continuity of $L$ on $\mathbb{R}_+^{\mathcal{R}}$ follows from Theorem 2.35, page 59 in Rockafellar and Wets [23] and the fact that it is convex, lower semicontinuous, and finite on $\mathbb{R}_+^{\mathcal{R}}$.

Under the stability condition (8), there exists some $\varepsilon > 0$ such that $(1+\varepsilon)\rho \in \mathcal{C}$. Thus,

$$\delta_K^*(x) \geq \sum_{r \in \mathcal{R}} x_r \log((1+\varepsilon)\rho_r).$$

It follows that

$$L(x) \geq \log(1+\varepsilon) \sum_{r \in \mathcal{R}} x_r.$$

This provides the first inequality in (9). In order to establish the second inequality, use the homogeneity property of $L$ to write

$$L(x) = \|x\|_\infty L(\|x\|_\infty^{-1} x) \leq \|x\|_\infty \sup_{y \in \mathbb{R}_+^{\mathcal{R}}, \|y\|_\infty = 1} L(y).$$

The supremum of a continuous function on a compact set is finite, which yields the second half of (9). □

It follows from equation (9) that, under condition (8), the stationary measure (7) can be normalized to a probability measure, which then implies stability (ergodicity) of the Markov process under the modified proportional fairness NBA, when (8) holds.

Fix now $y \in \mathbb{R}_+^{\mathcal{R}}$, and let $x = ny$, where $n$ is large. The heuristic calculation

$$\begin{aligned}\lambda_r^{\mathrm{PF}'}(x) &= \exp(n\delta_K^*(y) - n\delta_K^*(y - n^{-1}e_r)) \\ &\approx \exp(\partial_r \delta_K^*(y)) \\ &= \lambda_r^{\mathrm{PF}}(x),\end{aligned}$$

based on the homogeneity property of $\delta_K^*$, according to which $\delta_K^*(ny) = n\delta_K^*(y)$, and a heuristic Taylor approximation, suggests that the behavior of the systems under PF and PF' are similar, at least far from the origin. At this stage we content ourselves with making the following conjecture:

CONJECTURE 1. *Let $X^{\mathrm{PF}}$ and $X^{\mathrm{PF}'}$ denote the number of customers in steady state under PF and PF', respectively. We conjecture that the rescaled vectors $n^{-1}X^{\mathrm{PF}}$ and $n^{-1}X^{\mathrm{PF}'}$ satisfy, as $n \to \infty$, a large deviations principle with the same rate function $L$ as defined in (6).*

Remark that the vector of allocations $\lambda_r^{\mathrm{PF}'}(x)$ belongs to the convex set $\mathcal{C}$ for all $x \in \mathbb{Z}_+^{\mathcal{R}}$, in view of the following property of the function $\delta_K^*$:

LEMMA 4. *The function $\delta_K^*$ is such that, for all $x \in \mathbb{R}_+^{\mathcal{R}}$, and all $\varepsilon_r > 0$, $r \in \mathcal{R}$,*

$$\left\{ \frac{\delta_K^*(x) - \delta_K^*(x - \varepsilon_r e_r)}{\varepsilon_r} \right\} \in K. \qquad (10)$$

*It is understood in this expression that a vector $u$ with coordinates in $\{-\infty\} \cup \mathbb{R}$ belongs to $K$ when the vector $e^u$ with coordinates $e^{u_r}$ belongs to the original convex set $C$, and $e^{-\infty} = 0$.*

PROOF. Let $x \in \mathbb{R}_+^{\mathcal{R}}$. Assume first that all the coordinates $x_r$ are strictly positive. It follows that the vector $u$ achieving the supremum in the original definition of $\delta_K^*(x)$ is uniquely defined. By Lemma 1 above, and Theorem 25.1, page 242 in [22], it follows that $\delta_K^*$ is differentiable at $x$, its (ordinary) gradient being the vector $u$ achieving that supremum. Also, the function $\varepsilon \to \varepsilon^{-1}[\delta_K^*(x) - \delta_K^*(x - \varepsilon e_r)]$ is nonincreasing in $\varepsilon > 0$, and achieves its maximum as $\varepsilon \searrow 0$, where it equals the coordinate $u_r$ of the gradient (see Theorem 23.1, pages 213–214 in [22]). By monotonicity of the set $K$, it follows that $\delta_K^*$ satisfies the condition (10) at $x$.

We now show that the same is true when some coordinates of $x$ equal zero. Let $I \subset \mathcal{R}$ denote the set of indices $r$ for which the coordinate $x_r$ equals zero. We say that $x$ belongs to the *face $I$* when $x \in \mathbb{R}_+^{\mathcal{R}}$ and $x_r = 0$ if and only if $r \in I$. We also denote by $K_I$ the subset of $K$ consisting of these vectors $u$ with $u_r = -\infty$ if and only if $r \in I$. In the definition of $\delta_K^*(x)$, we may actually replace the optimization domain by $K_I$ rather than $K$. There is then a single vector $u$ of $K_I$ which achieves the corresponding supremum. We may conclude as in the previous case that (10) holds in the present case as well. $\square$

**3. Stability properties of proportional fairness.** The above characterization is now applied to the study of stability properties of the Markov process describing the number of users in the system under proportional fairness. Ergodicity is established by following the general approach of fluid limits, introduced in the contexts of more traditional queueing systems by Rybko and Stolyar [24] and Dai [9].

The section is organized as follows. The general model with Markovian users routing is first introduced. A characterization of the fluid limits of this process is then given. It is next established that the function $L$ defined in (6) is a Lyapunov function for these fluid limits, from which stability (or equivalently, ergodicity) of the original Markov process is deduced, under condition (8) for suitably defined loads $\rho_r$, $r \in \mathcal{R}$.

The model with Markovian users routing is as follows. As before, users are of different types, $r \in \mathcal{R}$. External arrivals of type $r$ users are according to a Poisson process with intensity $\overline{\nu}_r$; the service times of type-$r$ users are again exponential with parameter $\mu_r$. However, after completing service, type $r$ users will re-enter the system as type $s$ users with some probability $p_{rs}$. Thus the nonzero transition rates are now given by

$$
\begin{aligned}
x &\to x + e_r && \text{with rate } \overline{\nu}_r, \\
(11) \quad x &\to x - e_r + e_s && \text{with rate } \mu_r \lambda_r^{\mathrm{PF}}(x) p_{rs}, \\
x &\to x - e_r && \text{with rate } \mu_r \lambda_r^{\mathrm{PF}}(x) \left(1 - \sum_{s \in \mathcal{R}} p_{rs}\right).
\end{aligned}
$$

In the above, $e_r$ denotes the $r$th unit vector in $\mathbb{R}^{\mathcal{R}}$. It is assumed that the matrix $P = (p_{rs})_{r,s \in \mathcal{R}}$ is substochastic, and that its spectral radius is strictly less than 1. Thus, there exists a unique vector $\nu = (\nu_r)_{r \in \mathcal{R}}$ solving the traffic equations

$$\nu_r = \overline{\nu}_r + \sum_{s \in \mathcal{R}} p_{sr} \nu_s, \qquad r \in \mathcal{R},$$

also written in matrix form

$$(I - P^T)\nu = \overline{\nu},$$

where $P^T$ is the transposition of the routing probability matrix $P$. Introduce the notation $\rho_r = \nu_r / \mu_r$, and $\rho = (\rho_r)_{r \in \mathcal{R}}$. The main result of this section is the following:

THEOREM 1. *The Markov process with Markovian users routing is ergodic under condition* (8).

In order to establish the theorem, a characterization of the *fluid limits* of the original Markov process is required. To this end, the following definition will be used. Note that the constant $A$ appearing in this definition differs from the one appearing in Lemma 3. In the sequel, to simplify notations, $A$ will always be used to denote an arbitrary finite constant, whose value may vary from one statement to another.

DEFINITION 1. The functions $x_r : \mathbb{R}_+ \to \mathbb{R}_+$, $r \in \mathcal{R}$, are called fluid trajectories of the system with Markovian users routing if there exist non-decreasing, Lipschitz continuous functions $D_r : \mathbb{R}_+ \to \mathbb{R}_+$, $r \in \mathcal{R}$, such that $D_r(0) = 0$, admitting $A$ as a Lipschitz constant for any $A$ such that $\mathcal{C} \subset [0, A]^{\mathcal{R}}$, that verify

$$(12) \quad x_r(t) = x_r(0) + \overline{\nu}_r t - \mu_r D_r(t) + \sum_{s \in \mathcal{R}} p_{sr} \mu_s D_s(t), \qquad t \in \mathbb{R}_+, r \in \mathcal{R},$$

and for almost every $t \in \mathbb{R}_+$, all $r \in \mathcal{R}$, the derivatives $\dot{D}_r(t)$ exist and verify

$$(13) \qquad \dot{D}_r(t) \in \left[0, \limsup_{y \to x(t)} \lambda_r^{\mathrm{PF}}(y)\right],$$

$$(14) \qquad x_r(t) > 0 \Rightarrow \dot{D}_r(t) = \lambda_r^{\mathrm{PF}}(x(t)) = \exp(\gamma_r^{\mathrm{PF}}(x(t))).$$

The following notation will be used in the sequel. For any $x \in \mathbb{R}_+^{\mathcal{R}}$, $S(x)$ denotes the set of all fluid trajectories of the system with initial condition $x$. Thus it is a subset of $\mathcal{C}([0, +\infty), \mathbb{R}_+^{\mathcal{R}})$, that is the space of continuous, $\mathbb{R}_+^{\mathcal{R}}$-valued functions on $[0, +\infty)$.

Note that at this stage neither existence nor uniqueness of fluid trajectories with a given initial condition have been established.

The following result is the first step of the proof of Theorem 1. It implies as a corollary that the set $S(x)$ is nonempty, for any $x \in \mathbb{R}_+^{\mathcal{R}}$. However no claim of uniqueness of fluid trajectories is made.

THEOREM 2. *Consider a sequence of initial conditions* $X^k(0) = (X_r^k(0))_{r \in \mathcal{R}}$, $k \geq 1$, *such that for a sequence of positive numbers* $(z_k)_{k \in \mathbb{N}}$, $\lim_{k \to \infty} z_k = +\infty$, *and the limit* $\lim_{k \to \infty} z_k^{-1} X^k(0) = x(0)$ *exists in* $\mathbb{R}_+^{\mathcal{R}}$.

*Then for all $T > 0$, and all $\varepsilon > 0$, the following convergence takes place:*

$$\lim_{k \to \infty} \mathbf{P}\left( \inf_{f \in S(x(0))} \sup_{t \in [0,T]} |z_k^{-1} X^k(z_k t) - f(t)| \geq \varepsilon \right) = 0.$$

*In words, the restriction of the rescaled process $z_k^{-1} X^k(z_k \cdot)$ to any compact interval $[0, T]$ converges in probability to the set $S(x(0))$ of fluid trajectories* *with initial condition $x(0)$, where convergence of processes is for the uniform norm.*

The proof of Theorem 2 is deferred to the Appendix A. We expect a similar result to hold for other NBA, in particular for $\alpha$-fair NBA, provided one replaces $\lambda^{\text{PF}}$ by the corresponding allocation vector $\lambda^{\text{NBA}}$ in the definition of the fluid trajectories. Indeed the proof given in the Appendix A relies on two technical lemmas by Ye, Ou and Yuan [30] which apply to general $\alpha$-fair NBA, and the rest of the proof can be adapted in a straightforward manner.

The second step of the proof of Theorem 1 consists in establishing a suitable uniform convergence to zero of fluid trajectories:

THEOREM 3. *Under the stability condition* (8), *there exists $\tau > 0$ and $\varepsilon > 0$ such that, for any fluid trajectory $\{x(t)\}_{t \in \mathbb{R}_+}$, provided $L(x(0)) = 1$, then $L(x(\tau)) \leq 1 - \varepsilon$.*

The following lemma will be needed in the proof of Theorem 3:

LEMMA 5. *Let $\{x(t)\}_{t \in \mathbb{R}_+}$ be a fluid trajectory as per Definition 1. For every $t \geq 0$, let $I(t)$ denote the set of indices $r \in \mathcal{R}$ such that $x_r(t) = 0$, and $\bar{I}(t) = \mathcal{R} \setminus I(t)$.*

*(i) There exist modified arrival rates $\tilde{\nu}_r$, $r \in \bar{I}(t)$, and modified routing probabilities, $\tilde{p}_{rs}$, $r, s \in \bar{I}(t)$, that depend only on the set $I(t)$, such that the matrix $(\tilde{p}_{rs})_{r,s \in \bar{I}(t)}$ is sub-stochastic with spectral radius strictly less than 1, the identity*

$$(15) \qquad (\nu_r)_{r \in \bar{I}(t)} = (I - \tilde{P}^T)^{-1} \tilde{\nu}$$

*holds, and furthermore, for almost every $t > 0$,*

$$(16) \quad \begin{aligned} \frac{d}{dt} x_r(t) &= \tilde{\nu}_r + \sum_{r \in \bar{I}(t)} \mu_s \tilde{p}_{sr} \lambda_s^{\text{PF}}(x(t)) - \mu_r \lambda_r^{\text{PF}}(x(t)), & r \in \bar{I}(t), \\ \frac{d}{dt} x(t) &= 0, & r \in I(t). \end{aligned}$$

*Let $f(t) := L(x(t))$.*

*(ii) For almost every $t > 0$, it holds that:*

$$(17) \qquad \limsup_{h \searrow 0} \frac{f(t+h) - f(t)}{h} \leq \sum_{r \in \bar{I}(t)} (\gamma_r^{\text{PF}}(x(t)) - \log(\rho_r)) \dot{x}_r(t),$$

*where the derivatives $\dot{x}_r(t)$ are as in* (16).

(iii) *There exists a constant $A$ such that, for all $t > 0$,*

(18) $$\limsup_{h \searrow 0} \frac{f(t+h) - f(t)}{h} \leq A.$$

The proof of the lemma is given in Appendix B.
The following auxiliary result will also be used:

LEMMA 6. *Let a continuous function $f:[0,T] \to \mathbb{R}$ be given. Assume that there exists $\varepsilon \in \mathbb{R}$ such that, for almost all $t \in [0,T]$:*

(19) $$\limsup_{h \searrow 0} \frac{f(t+h) - f(t)}{h} \leq -\varepsilon.$$

*Assume further the existence of a constant $A \in \mathbb{R}$ such that for all $t \in [0,T]$,*

(20) $$\limsup_{h \searrow 0} \frac{f(t+h) - f(t)}{h} \leq A.$$

*Then it holds that, for all $s, t \in [0, T)$, $s < t$,*

(21) $$f(t) - f(s) \leq -\varepsilon(t - s).$$

REMARK 1. The following example illustrates the role of assumption (20) in Lemma 6. Let $f^+(t) = m([0,t])$, where $m$ is the uniform measure on the Cantor set obtained by successive exclusion of the middle third from the interval $[0, 1]$ (see, e.g., Falconer [11] for background). More precisely, this measure can be defined by specifying the mass it puts on intervals $[0, x]$ where $x$ is a *triadic* number, that is

$$x = \sum_{i=1}^{\infty} z_i 3^{-i},$$

where $z_i \in \{0, 1, 2\}$, $i \geq 1$. The uniform measure $m$ on this Cantor set is then specified by

$$m([0, x]) = \sum_{i=1}^{k} z_i 2^{-i-1},$$

where $k = \min\{i \geq 1 : z_i = 1\}$.

Define then

$$f(t) = -\varepsilon t + f^+(t).$$

The function $f$ is continuous, because the measure $m$ has no atoms. Moreover, the measure $m$ is supported by a set of null Lebesgue measure, so that the function $f$ satisfies condition (19) of Lemma 6. However, the conclusion (21) does not hold, precisely because condition (20) is not satisfied.

The result of Theorem 3 is established as follows.

PROOF OF THEOREM 3. Let $\{x(t)\}_{t\in\mathbb{R}_+}$ denote a fluid trajectory. Introduce the notation $u_r = \log(\lambda_r^{\mathrm{PF}}(x(t))/\rho_r)$. The right-hand side of equation (17), which we shall denote $h(t)$, can then be rewritten, in view of (16), as

$$h(t) = \sum_{r\in\bar{I}(t)} u_r \left[\tilde{\nu}_r - \nu_r e^{u_r} + \sum_{s\in\bar{I}(t)} \tilde{p}_{sr}\nu_s e^{u_s}\right],$$

or equivalently, in matrix form,

$$h(t) = \langle u, \tilde{\nu} - (I - \tilde{P}^T)(\nu e^u)\rangle.$$

We use the notation $\lceil\eta\rfloor$ to denote the diagonal matrix with diagonal entries provided by the coordinates of the vector $\eta$. Elementary manipulations entail that

$$\begin{aligned}
h(t) &= -\langle u, (I - \tilde{P}^T)\lceil\nu\rfloor(e^u - 1)\rangle \\
&= -\langle \nu, \lceil(I - \tilde{P})u\rfloor(e^u - 1)\rangle \\
&= -\langle \tilde{\nu}, (I - \tilde{P})^{-1}\lceil(I - \tilde{P})u\rfloor(e^u - 1)\rangle,
\end{aligned} \quad (22)$$

the first equality relying on identity (15). In order to show that the previous expression is nonpositive, it is enough to show that for each $r \in \bar{I}(t)$, the coefficient of $\tilde{\nu}_r$ is nonpositive, that is,

$$(23) \quad F_r(u) := \sum_{n\geq 0}\sum_{s\in\bar{I}(t)} \tilde{p}_{rs}^{(n)}(e^{u_s} - 1)\left[u_s - \sum_{\ell\in\bar{I}(t)} p_{s\ell}u_\ell\right] \geq 0, \qquad r \in \bar{I}(t).$$

The following lemma, whose proof is deferred to the Appendix D, is now needed:

LEMMA 7. *For any substochastic matrix $\tilde{P} = (\tilde{p}_{rs})_{r,s\in\bar{I}}$ with spectral radius strictly less than 1, and any real numbers $u_s$, $s \in \bar{I}$, then:*

(i) *Inequality* (23) *holds.*

(ii) *The function $F_r$ as defined in* (23) *verifies $F_r(u) \geq F_r(u^+)$ for all $u \in \mathbb{R}_+^{\bar{I}}$, where $u^+ := (u_s^+)_{s\in\bar{I}}$, and $u_s^+ = \max(u_s, 0)$.*

(iii) *There is equality in* (23) *only if for all states $s$ such that $\sum_{n\geq 0}\tilde{p}_{rs}^{(n)} > 0$, one has $u_s = 0$.*

That the term $h(t)$ is nonpositive follows from Lemma 7(i).

When $x(t) \neq 0$, the allocation vector $\lambda^{\mathrm{PF}}(x(t))$ must lie on the external boundary of the capacity set $\mathcal{C}$. Thus, by (8), for some positive $\varepsilon$, there must exist some coordinate $r$ such that $\lambda_r^{\mathrm{PF}}(x(t)) \geq (1+\varepsilon)\rho_r$. Therefore, setting $\delta = \log(1+\varepsilon) > 0$, it holds that $u_s \geq \delta$ for some $s \in \bar{I}$. There must also exist some $r \in \bar{I}$ such that $\tilde{\nu}_r > 0$, and $\sum_{n\geq 0} \tilde{p}_{rs}^{(n)} > 0$. It is also the case that the $u_k$ are bounded from above by some constant $A$, since the allocations $\lambda_k^{\mathrm{PF}}$ are bounded from above.

By Lemma 7(ii), one thus has

$$x(t) \neq 0 \Rightarrow h(t) \leq - \inf_{r:\tilde{\nu}_r>0} \tilde{\nu}_r \inf_{u \in S} F_r(u),$$

where the set $S$ is defined as

$$S = \left\{ u \in [0, A]^{\bar{I}} : \max_{k \in \bar{I}} u_k \geq \delta \right\}.$$

Since the function $F_r$ is continuous and the set $S$ is compact, the infimum of $F_r(u)$ over $S$ is attained; however it cannot be zero, in view of Lemma 7(iii) and the definition of $S$. Thus, the right-hand side of the above is less than $-\varepsilon(I(t))$ for some strictly positive $\varepsilon(I(t))$ that depends only on the set $I(t)$.

By assumption, the initial condition of the fluid trajectory in the statement of Theorem 2 is such that $L(x(0)) = 1$.

Thus, in view of (9), there exists $r$ so that $x_r(0) \geq 1/K$ for some finite positive constant $K$. Setting $\tau = 1/(2KA)$, where $A$ is such that the capacity set $\mathcal{C}$ is a subset of $[0, A]^{\mathcal{R}}$, it then follows that for any fluid trajectory with initial condition $x(0)$ such that $L(x(0)) = 1$, then $x(t) \neq 0$ on $[0, \tau]$. Hence, by the previous evaluations, in view of (17,18) for any such fluid trajectory, the function $f(t) = L(x(t))$ satisfies the assumptions of Lemma 6, with $\varepsilon := \inf_{I \subset \mathcal{R}, I \neq \mathcal{R}} \varepsilon(I) > 0$. Thus, by Lemma 6:

$$L(x(\tau)) \leq 1 - \tau\varepsilon < 1.$$

The claim of the theorem follows. $\square$

The proof of Theorem 1 will require to combine Theorems 2, 3 and the following ergodicity criterion, which is a direct consequence of Theorem 8.13, page 224 in Robert [21]:

THEOREM 4 ([21]). *Let $X(t)$ be a Markov jump process on a countable state space $\mathcal{S}$. Assume there exists a function $L : \mathcal{S} \to \mathbb{R}_+$ and constants $A$, $\varepsilon$, and an integrable stopping time $\hat{\tau} > 0$ such that for all $x \in \mathcal{S}$:*

(24) $$L(x) > A \Rightarrow \mathbf{E}_x L(X(\hat{\tau})) \leq L(x) - \varepsilon \mathbf{E}_x(\hat{\tau}).$$

*If in addition the set $\{x : L(x) \leq A\}$ is finite, and $\mathbf{E}_x L(X(1)) < +\infty$ for all $x \in \mathcal{S}$, then the process $X(t)$ is ergodic.*

PROOF OF THEOREM 1. Let $\tau$ be as in Theorem 3. Consider the deterministic stopping time $\hat{\tau} = L(x(0))\tau$. Denote by $\mathbf{P}_x$ the probability distribution of the Markov process $(X(t))$ with initial condition $x \in \mathbb{N}^{\mathcal{R}}$.

It is readily seen that the collection of probability distributions

$$\{\mathbf{P}_x(L(x)^{-1}X_r(\hat{\tau}) \in \cdot)\}_{r \in \mathcal{R}, x \in \mathbb{N}^{\mathcal{R}} \setminus \{0\}}$$

is uniformly integrable. Indeed, let $A_r$ denote independent unit rate Poisson processes, used to generate users arrival times. Then the process $X(t)$ can be generated so that

(25) $$X_r(t) \leq X_r(0) + A_r(\nu_r t), \qquad t \geq 0, r \in \mathcal{R}.$$

Thus, for $X(0) = x$,

$$\frac{X_r(\hat{\tau})}{L(x)} \leq \frac{x_r}{L(x)} + \frac{A_r(\nu_r L(x))}{L(x)}.$$

The first term is bounded from above uniformly in $x \neq 0$, in view of Lemma 3, (9). The second term has mean 1. Its variance equals $\nu_r/L(x)$. Thus the second moments of these variables are uniformly bounded in $x \neq 0$. Therefore, Lavallée–Poussin criterion for uniform integrability applies.

In view of (9), it then follows that the collection of probability distributions

$$\{\mathbf{P}_x(L(x)^{-1}L(X(\hat{\tau})) \in \cdot)\}_{x \in \mathbb{N}^{\mathcal{R}} \setminus \{0\}}$$

is also uniformly integrable.

The result of Theorem 2 entails that for any sequence of initial conditions $x^k$ such that $\|x^k\|_\infty \to \infty$ as $k \to \infty$, the corresponding rescaled variables $X^k(L(x^k)\tau)/L(x^k)$ converge in probability to the set $V$ defined as

$$V := \bigcup_{x \in \mathbb{R}_+^{\mathcal{R}}, L(x)=1} \{x(\tau), x(\cdot) \in S(x)\}.$$

In words, $V$ is the set of states of fluid trajectories at time $\tau$ for all fluid trajectories with initial condition $x(0)$ satisfying $L(x(0)) = 1$.

It can be verified from (9) and the definition of fluid trajectories that the set $V$ is compact. Continuity of $L$ together with compactness of $V$ entail that the sequence of random variables $L(X^k(L(x^k)\tau))/L(x^k)$ converges in probability to the set $L(V)$.

Thus, by Theorem 3, the sequence of random variables $L(X^k(L(x^k)\tau))/L(x^k)$ converges in probability to the interval $[0, 1-\varepsilon]$, where $\varepsilon > 0$.

Combined with the uniform integrability just shown, this yields

$$\limsup_{L(x) \to \infty} \frac{1}{L(x)} \mathbf{E}_x L(X(L(x)\tau)) \leq 1 - \varepsilon.$$

Thus, Condition (24) of Theorem 4 holds for $A$ sufficiently large. The second requirement, that the set $\{x : L(x) \leq A\}$ be finite, follows from (9). Finally, the last condition, that is $\mathbf{E}_x L(X(1)) < +\infty$ for all $x$ is easily verified, invoking once more the bounds (25) and (9). $\square$

**4. Application to phase-type service distributions.** We now apply Theorem 1 to systems with general phase-type distributions rather than exponential service time distributions. More precisely, we consider the same setting as before, with user classes $r \in \mathcal{R}$, and capacity set $\mathcal{C} \subset \mathbb{R}_+^{\mathcal{R}}$. New type $r$ users arrive as usual according to a Poisson process with intensity $\nu_r$.

The service time distribution of type $r$ customers is now defined as follows. A finite set $I_r$, referred to as the set of service phases, is given. The total service time is characterized as the aggregation of service times required in subsequent visits to phases. At each visit to phase $i$, a corresponding service time that is exponentially distributed, with parameter $\mu_{r,i}$, is required. A visit to phase $i$ is followed by a visit to phase $j$ with probability $p_{r;ij}$. A probability distribution $\{\alpha_i\}_{i \in I_r}$ on $I_r$ specifies the phase in which service starts. The transition matrix $P_r := (p_{r;ij})_{i,j \in I_r}$ is assumed to be sub-stochastic, with spectral radius strictly less than 1. It is easily checked that the above description is equivalent to the definition of phase-type distributions given in [1], page 83.

Denote by $\widehat{\mathcal{R}}$ the set of pairs $(r,i)$ with $r \in \mathcal{R}$ and $i \in I_r$. For all $(r,i) \in \widehat{\mathcal{R}}$, let $x_{r,i}$ denote the number of class $r$ users who are currently in phase $i$ of their service.

The process $(x_{r,i})_{(r,i) \in \widehat{\mathcal{R}}}$ is then a Markov process of the kind covered by Theorem 1. More precisely, it corresponds to the following parameters. For the class $s = (r,i) \in \widehat{\mathcal{R}}$, the external arrival rate $\overline{\nu}_s$ is given by $\nu_r \alpha_{r,i}$ and the corresponding service time parameter is $\mu_s = \mu_{r,i}$. For two classes $s = (r,i)$, $s' = (r',i')$, the corresponding routing probability $p_{ss'}$ is zero if $r \neq r'$, and otherwise equals $p_{r;ii'}$. Finally, the capacity set $\widehat{\mathcal{C}}$ is determined from the original capacity set $\mathcal{C}$ as follows. The allocation vector $(\lambda_s)_{s \in \widehat{\mathcal{R}}}$ belongs to $\widehat{\mathcal{C}}$ if and only if the allocation vector $(\lambda_r)_{r \in \mathcal{R}}$ belongs to $\mathcal{C}$, where $\lambda_r$ is given by $\sum_{i \in I_r} \lambda_{(r,i)}$.

We then have the following:

THEOREM 5. *The process tracking the numbers $x_r$ of users of class $r$, under proportionally fair allocation of resources characterized by the set $\mathcal{C}$, assuming Poisson arrivals and phase type distributions as just described, is ergodic under the usual condition* (8), *where $\rho_r = \nu_r \sigma_r$, and $\sigma_r$ is the mean service time for class $r$ users.*

PROOF. By Theorem 1, ergodicity holds provided the vector $(\rho_{(r,i)})_{(r,i)\in\widehat{\mathcal{R}}}$ belongs to the interior of $\hat{\mathcal{C}}$. Equivalently, it holds if the vector with $r$th coordinate $\sum_{i\in I_r}\rho_{(r,i)}$ belongs to $\overset{\circ}{\mathcal{C}}$.

With the specific routing probability matrix $P$ obtained from the characteristics of the phase type service distributions, one has

$$\rho_{(r,i)} = \frac{1}{\mu_{r,i}} \sum_{j\in I_r} \sum_{n\geq 0} p^{(n)}_{r;ji}\overline{\nu}_{r,j}$$

$$= \frac{1}{\mu_{r,i}} \sum_{j\in I_r} \sum_{n\geq 0} p^{(n)}_{r;ji}\alpha_{r,j}\nu_r.$$

This in turn implies that

$$\sum_{i\in I_r}\rho_{(r,i)} = \nu_r \sum_{i\in I_r} \frac{1}{\mu_{r,i}} \sum_{j\in I_r, n\geq 0} \alpha_{r,j}\, p^{(n)}_{r;ji}.$$

Noting that in the above expression, the last sum over $j\in I_r$ and $n\geq 0$ gives the average number of visits to phase $i$ in a class $r$ service time, it readily follows that this last expression coincides with $\nu_r\sigma_r$, which completes the proof. $\square$

**5. Relationships between balanced fairness and proportional fairness.** In this section we define the balanced fairness NBA, give an equivalent characterization and then use it to relate the stationary distributions under balanced fairness and modified proportional fairness.

The balanced fairness NBA, introduced in [5], is best defined in terms of the *balance function*. The balance function, denoted $\psi$, is defined by induction on $\mathbb{Z}_+^{\mathcal{R}}$, starting from $\psi(0)=1$, $\psi(x)=0$ for any $x$ not in $\mathbb{R}_+^{\mathcal{R}}$, and

$$\psi(x) = \inf\{a > 0 : \{a^{-1}\psi(x-e_r)\}_{r\in\mathcal{R}} \in C\},$$

where $e_r$ is the $r$th unit vector in $\mathbb{R}^{\mathcal{R}}$. The balanced fairness rate allocation vector $\lambda^{\mathrm{BF}}$ is then defined as

$$\lambda_r^{\mathrm{BF}}(x) = \frac{\psi(x-e_r)}{\psi(x)}, \qquad x\in\mathbb{Z}_+^{\mathcal{R}}, r\in\mathcal{R}.$$

As for proportional fairness, it is convenient to consider the logarithms $\gamma_r$ of the allocated capacities $\lambda_r$, rather than the $\lambda_r$ themselves. Denote by $\gamma^{\mathrm{BF}}(x)$ the vector of logarithms of balanced fair allocations, that must lie in the convex nonincreasing set $K$. Introduce the notation $\phi(x) = -\log\psi(x)$. Thus one has

$$\gamma_r^{\mathrm{BF}}(x) = \phi(x) - \phi(x - e_r).$$

In other words, the vector $\gamma^{\mathrm{BF}}(x)$ is given by the increments of the function $\phi$ at $x$, and can be seen as an approximate gradient of $\phi$ at $x$. We introduce the notation $\nabla_d f(x) = \{f(x) - f(x - e_r)\}_{r \in \mathcal{R}}$. Note that a stationary measure for the Markov process counting users of all types is given in terms of the function $\phi$ by

$$\pi^{\mathrm{BF}}(x) = \frac{1}{Z} \exp\left(-\phi(x) + \sum_{r \in \mathcal{R}} x_r \log(\rho_r)\right) \tag{26}$$

for some normalization constant $Z$. A consequence of the reversibility property of the Markov process is that this measure is also stationary for the modified Markov process with Markovian routing [5].

We now give an alternative definition of $\phi$.

LEMMA 8. *The function $\phi$ admits the following characterization:*

$$\phi(x) = \sup_{f \in \mathcal{F}} \{f(x)\}$$

*where $\mathcal{F}$ is the set of functions defined on $\mathbb{Z}^{\mathcal{R}}$ such that $f(0) = 0$, $f(y) = +\infty$ for $y \notin \mathbb{Z}_+^{\mathcal{R}}$, and $\nabla_d f(y)$ belongs to $K$ for all $y \in \mathbb{Z}_+^{\mathcal{R}}$.*

PROOF. Denote by $\hat{\phi}(x)$ the result of the optimization problem in the right-hand side of the above expression. Proceed by induction on $x \in \mathbb{Z}_+^{\mathcal{R}}$ to show that $\phi(x) = \hat{\phi}(x)$. Clearly, $\phi(0) = \hat{\phi}(0) = 0$. Also, as the function $\phi$ satisfies the conditions over which the optimization is performed, necessarily one has that $\phi(x) \leq \hat{\phi}(x)$, for all $x \in \mathbb{Z}_+^{\mathcal{R}}$. Assume thus that $\hat{\phi}(y) = \phi(y)$ for all $y \leq x$, $y \neq x$. The definition by induction of $\psi$ implies that

$$\phi(x) = \sup\{a : \{a - \phi(x - e_r)\}_{r \in \mathcal{R}} \in K\}.$$

On the other hand, for any $f$ satisfying the assumptions,

$$f(x) \leq \sup\{a : \{a - f(x - e_r)\}_{r \in \mathcal{R}} \in K\}$$
$$\leq \sup\{a : \{a - \phi(x - e_r)\}_{r \in \mathcal{R}} \in K\}$$
$$= \phi(x).$$

We have used for the first inequality the definition of the constraints satisfied by $f$, for the second we have used the induction hypothesis that $f(x - e_r) \leq \phi(x - e_r)$ together with monotonicity of the set $K$, and the last equality is just the inductive definition of $\phi$. □

We are now ready to establish the following:

THEOREM 6. *For any $x \in \mathbb{Z}_+^{\mathcal{R}}$, the following inequalities hold:*

$$\delta_K^*(x) \leq \phi(x) \leq \delta_K^*(x) + r(x), \tag{27}$$

*where*

$$r(x) := \sum_{r \in \mathcal{R}: x_r > 0} \sum_{m=1}^{x_r} \frac{1}{m}.$$

PROOF. The two inequalities shall be established by induction on $\sum_{r \in \mathcal{R}} x_r$ for $x \in \mathbb{Z}_+^{\mathcal{R}}$. They obviously hold true for $x = 0$, as $\phi(0) = \delta_K^*(0) = 0$. Assume thus that they hold for all $y \in \mathbb{Z}_+^{\mathcal{R}}$ such that $\sum_{r \in \mathcal{R}} y_r \leq n$, for some integer $n \geq 0$, and let $x \in \mathbb{Z}_+^{\mathcal{R}}$ be given, $\sum_{r \in \mathcal{R}} x_r = n+1$. By the induction hypothesis and the result of Lemma 8, it holds that

$$\phi(x) = \sup\{a : \{a - \phi(x - e_r)\}_{r \in \mathcal{R}} \in K\}$$
$$\geq \sup\{a : \{a - \delta_K^*(x - e_r)\}_{r \in \mathcal{R}} \in K\}.$$

Now, in view of Lemma 4, it holds that

$$\{\delta_K^*(x) - \delta_K^*(x - e_r)\}_{r \in \mathcal{R}} \in K.$$

Therefore,

$$\phi(x) \geq \delta_K^*(x),$$

and the first inequality in (27) is established.

By the induction hypothesis again, we have that

$$\phi(x) \leq \sup\{a : \{a - \delta_K^*(x - e_s) - r(x - e_s)\}_{s \in \mathcal{R}} \in K\}. \tag{28}$$

Consider first the case where $x_s > 0$ for all $s \in \mathcal{R}$. We shall rely on the following lemma, the proof of which will be given after the end of the current proof.

LEMMA 9. *For all $x, h \in \mathbb{R}^{\mathcal{R}}$, such that $x$ has strictly positive coordinates, and $x + h$ has nonnegative coordinates, it holds that*

$$\delta_K^*(x + h) \leq \delta_K^*(x) + \langle h, \gamma^{\mathrm{PF}}(x) \rangle + \sum_{s \in \mathcal{R}} \frac{h_s^2}{x_s}. \tag{29}$$

Thus, in view of the previous equation, we have that

$$\delta_K^*(x - e_s) \leq \delta_K^*(x) - \gamma_s^{\mathrm{PF}}(x) + \frac{1}{x_s}.$$

Combining this upper bound with (28), as the vector $(\gamma_s^{\mathrm{PF}}(x))_{s \in \mathcal{R}}$ is in $K$, we have that

$$\phi(x) \leq \delta_K^*(x) + \sup_{s \in \mathcal{R}} \left\{ r(x - e_s) + \frac{1}{x_s} \right\}.$$

In view of the definition of $r(y)$, the second term in the right-hand side is clearly upper bounded by $r(x)$, which establishes the desired inequality for $x$.

To conclude the proof, it remains to deal with the case where some coordinates $x_s$ equal zero. This case is in fact similar to the previous one: if $x$ belongs to face $I$ (i.e., $x_s = 0$ if and only if $s \in I$), the previous argument carries over in $\mathbb{Z}_+^{\mathcal{R}\setminus I}$ by considering the convex set $K_I$ instead of $K$. □

PROOF OF LEMMA 9. Let $x, h \in \mathbb{R}^{\mathcal{R}}$ be fixed, such that $x$ has strictly positive coordinates, and $x + h$ has nonnegative coordinates. Let $\gamma \in \mathbb{R}^{\mathcal{R}}$ be such that
$$\delta_K^*(x) = \langle x, \gamma \rangle - \delta_K(\gamma).$$
The pair $(x, \gamma)$ verifies the relations $x \in \partial \delta_K(\gamma)$, $\gamma \in \partial \delta_K^*(x)$. In addition, the following one-to-one correspondence between subgradients of $\delta_{\mathcal{C}}$ and $\delta_K$ can be established:
$$x \in \partial \delta_K(\gamma) \Leftrightarrow \{x_s e^{-\gamma_s}\}_{s \in \mathcal{R}} \in \partial \delta_{\mathcal{C}}(\{e^{\gamma_s}\}_{s \in \mathcal{R}}).$$
Let $h \in \mathbb{R}^{\mathcal{R}}$ be fixed. We have that
$$\delta_K^*(x + h) = \sup_{g \in \mathbb{R}^{\mathcal{R}}} \{\langle g, x + h \rangle - \delta_K(g)\}$$
$$= \sup_{u \in \mathbb{R}^{\mathcal{R}}} \{\langle u + \gamma, x + h \rangle - \delta_K(u + \gamma)\}$$
$$= \delta_K^*(x) + \langle h, \gamma \rangle + \sup_{u \in \mathbb{R}^{\mathcal{R}}} \{\langle u, x + h \rangle + \delta_K(\gamma) - \delta_K(\gamma + u)\}.$$
However, by convexity of $\delta_{\mathcal{C}}$, and recalling that $e^{-\gamma} x \in \partial \delta_{\mathcal{C}}(e^\gamma)$, we have the following inequality:
$$\delta_K(\gamma + u) = \delta_{\mathcal{C}}(e^{\gamma+u})$$
$$\geq \delta_{\mathcal{C}}(e^\gamma) + \langle e^{\gamma+u} - e^\gamma, e^{-\gamma} x \rangle$$
$$= \delta_K(\gamma) + \sum_{s \in \mathcal{R}} (e^{u_s} - 1) x_s.$$
Combined with the previous expression for $\delta_K^*(x + h)$, this yields
$$\delta_K^*(x + h) \leq \delta_K^*(x) + \langle h, \gamma \rangle + \sup_{u \in \mathbb{R}^{\mathcal{R}}} \left\{ \sum_{s \in \mathcal{R}} u_s(x_s + h_s) - (e^{u_s} - 1) x_s \right\}$$
$$= \delta_K^*(x) + \langle h, \gamma \rangle + \sum_{s \in \mathcal{R}: x_s + h_s > 0} (x_s + h_s) \log(1 + h_s/x_s) - h_s$$
$$+ \sum_{s \in \mathcal{R}: x_s + h_s = 0} x_s.$$

The claimed inequality (29) now follows by noting that (i) $\log(1+h_s/x_s) \leq h_s/x_s$, and (ii) $\gamma = \gamma^{\mathrm{PF}}(x)$. □

A simple consequence of the theorem is the following.

COROLLARY 1. *For any $x \in \mathbb{R}^+$, it holds that*
$$\lim_{n \to \infty} \frac{1}{n} \phi(nx) = \delta_K^*(x).$$

PROOF. This follows trivially since the function $\delta_K^*$ is positively homogeneous, that is, $\delta_K^*(nx) = n \delta_K^*(x)$, and since the remainder term $r(nx)$ in (27) is of order $\log(n)$, and a fortiori $o(n)$. □

Note that, in view of (26),
$$\frac{1}{n} \log(\pi^{\mathrm{BF}}(nx)) = -\frac{\phi(nx)}{n} + \sum_{r \in \mathcal{R}} x_r \log(\rho_r) - \frac{\log(Z)}{n}.$$

The last term must go to zero as $n \to \infty$. This together with Corollary 1 yield the following:

COROLLARY 2. *The stationary distribution $\pi^{\mathrm{BF}}$ as in (26) admits the following large deviations asymptotics:*
$$(30) \qquad \lim_{n \to \infty} \frac{1}{n} \log \pi^{\mathrm{BF}}(nx) = -L(x), \qquad x \in \mathbb{R}_+^{\mathcal{R}},$$

*where $L$ is the Lyapunov function (6) used in the study of stability properties of proportional fairness. It thus admits the same large deviations characteristics as the stationary distribution (7) of the system under PF' sharing.*

REMARK 2. The result of Theorem 6 also implies that, if for all $x \in \mathbb{R}_+^{\mathcal{R}}$, there exists a limit $\lim_{n \to \infty} \lambda^{\mathrm{BF}}(nx)$ of the allocation vector under balanced fairness, then it must coincide with $\lambda^{\mathrm{PF}}(x)$. So far we have not been able to establish the existence of such a limit, except in the special case where $|\mathcal{R}| = 2$, although it seems plausible that the limit exists more generally.

## APPENDIX A: PROOF OF THEOREM 2

We argue by contradiction, assuming that for some $\varepsilon > 0$, and for all $k$ in an infinite subsequence of the original sequence, it holds that

$$(31) \qquad \mathbf{P}\left(\inf_{f \in S(x(0))} \sup_{t \in [0,T]} |z_k^{-1} X^k(z_k t) - f(t)| \geq \varepsilon \right) \geq \varepsilon.$$

In the rest of the proof, without loss of generality we assume that the above evaluation holds true for all $k \geq 1$.

The trajectories of the processes $X^k(t)$ can be represented explicitly in terms of independent unit rate Poisson processes $A_r^k$, $D_{rs}^k$, $D_r^k$, $k \geq 0$, $r$, $s \in \mathcal{R}$, as follows:

$$X_r^k(t) = X_r^k(0) + A_r^k(\overline{\nu}_r t) + \sum_{s \in \mathcal{R}} D_{sr}^k\left(\mu_s p_{sr} \int_0^t \lambda_s^{\mathrm{PF}}(X^k(u))\,du\right)$$

$$- \sum_{s \in \mathcal{R}} D_{rs}^k\left(\mu_r p_{rs} \int_0^t \lambda_r^{\mathrm{PF}}(X^k(u))\,du\right)$$

$$- D_r^k\left(\mu_r\left(1 - \sum_{s \in \mathcal{R}} p_{rs}\right) \int_0^t \lambda_r^{\mathrm{PF}}(X^k(u))\,du\right).$$

This implies the following, by a change of variables in the integrals, and using the fact that $\lambda^{\mathrm{PF}}(ax) = \lambda^{\mathrm{PF}}(x)$ for all scalar $a > 0$:

$$\frac{1}{z_k} X_r^k(z_k t) = \frac{1}{z_k} X_r^k(0) + \overline{\nu}_r t + \sum_{s \in \mathcal{R}} \mu_s p_{sr} \int_0^t \lambda_s^{\mathrm{PF}}(z_k^{-1} X^k(z_k u))\,du$$

$$- \mu_r \int_0^t \lambda_r^{\mathrm{PF}}(z_k^{-1} X^k(z_k u))\,du + \varepsilon_r^k(t),$$

where the error term $\varepsilon_r^k(t)$ verifies, for all $T > 0$,

$$\sup_{t \in [0,T]} |\varepsilon_r^k(t)| \leq \frac{1}{z_k} \sup_{t \in [0,\overline{\nu}_r T]} |A_r^k(z_k t) - z_k t|$$

$$+ \sum_{s \in \mathcal{R}} \frac{1}{z_k} \sup_{t \in [0,\mu_s p_{sr} AT]} |D_{sr}^k(z_k t) - z_k t|$$

$$+ \sum_{s \in \mathcal{R}} \frac{1}{z_k} \sup_{t \in [0,\mu_r p_{rs} AT]} |D_{rs}^k(z_k t) - z_k t|$$

$$+ \frac{1}{z_k} \sup_{t \in [0,\mu_r AT]} |D_r^k(z_k t) - z_k t|.$$

In these expressions, $A$ is a constant such that $\mathcal{C} \subset [0,A]^{\mathcal{R}}$.

The following large deviations bound on the maximal deviation of a unit rate Poisson process from its mean is now needed:

LEMMA A.1. *Let $\Xi$ be a unit rate Poisson process. Then for all $T > 0$, and all $\lambda > 0$, it holds that*

$$(32) \qquad \mathbf{P}\left(\sup_{0 \leq t \leq T} |\Xi(t) - t| \geq \lambda T\right) \leq e^{-Th(\lambda)} + e^{-Th(-\lambda)},$$

*where*

$$(33) \qquad h(\lambda) := (1+\lambda)\log(1+\lambda) - \lambda$$

is the Cramér transform of a unit mean, centered Poisson random variable. In the above formula, it is understood that $h(-\lambda) = +\infty$ if $\lambda > 1$.

This result and its proof are standard (see, e.g., [26]). It implies the existence of a subsequence $k(\ell), \ell \geq 1$ of the original sequence, and a sequence $\varepsilon(\ell)$ decreasing to zero, such that

$$\sum_{\ell \geq 1} \mathbf{P}\bigg(\sup_{t \in [0,T]} |\varepsilon_r^{k(\ell)}(t)| \geq \varepsilon(\ell)\bigg) < \infty, \qquad r \in \mathcal{R}.$$

Indeed, it can be checked from the definition of $\varepsilon_r^k(t)$ and Lemma A.1 that the sum in the left-hand side is finite for the particular choice

$$\begin{cases} k(1) = 1, \\ k(\ell) = \min\{k > k(\ell-1) : z_k \geq \ell\}, & \ell > 1, \\ \varepsilon(\ell) = \ell^{-1/4}, & \ell \geq 1. \end{cases}$$

Without loss of generality we again assume that finiteness of the sum holds true for the original sequence $k \geq 1$. Thus, by Borel–Cantelli's lemma, almost surely, $\sup_{t \in [0,T]} |\varepsilon_r^k(t)| \to 0$ as $k \to \infty$. The following variation on Arzela–Ascoli's theorem will then be used to proceed:

LEMMA A.2 (Lemma 6.3, [30]). *Suppose that a sequence of functions $f_k : [0,T] \to \mathbb{R}$ has the following properties:*

(i) $\{f_k(0)\}_{k \geq 0}$ *is bounded;*
(ii) *there is a constant $M > 0$, and a sequence of positive numbers $\sigma_k$, with $\sigma_k \to 0$ as $k \to \infty$, such that*

$$|f_k(t) - f_k(s)| \leq M(t-s) + \sigma_k, \qquad k \geq 0, s, t \in [0,T].$$

*Then the sequence admits a subsequence that converges uniformly on $[0,T]$ to a Lipschitz continuous function $f : [0,T] \to \mathbb{R}$ with Lipschitz constant $M$.*

In the present setup, this result guarantees that for any $T > 0$, with probability 1, for any subsequence of the original sequence $k \geq 1$, there exists a further subsequence, denoted $k'$, along which, for all $r \in \mathcal{R}$, the following convergences take place, uniformly on $[0,T]$:

$$\int_0^t \lambda_r^{\mathrm{PF}}(z_{k'}^{-1} X^{k'}(z_{k'} u)) \, du \to D_r(t),$$

$$z_{k'}^{-1} X_r^{k'}(z_{k'} t) \to x_r(t) := x(0) + \overline{\nu}_r t + \sum_{s \in \mathcal{R}} p_{sr} \mu_s D_s(t) - \mu_r D_r(t),$$

where the functions $D_r$ are $A$-Lipschitz. (A set in which any infinite sequence admits a convergent subsequence is usually called *sequentially compact.*

Sequential compactness is equivalent to compactness in the case of metric spaces.) We shall now establish that all such limits are fluid trajectories of the system. To this end, the following lemma, also taken from [30], Lemma 6.2(b), is needed:

LEMMA A.3. *For all $r \in \mathcal{R}$, and any $x \in \mathbb{R}_+^{\mathcal{R}}$ such that $x_r > 0$, the bandwidth allocation function $\lambda_r^{\mathrm{PF}}$ is continuous at $x$.*

In fact, Ye, Ou and Yuan establish this result in the context of particular, polyhedral capacity sets $\mathcal{C}$; however their proof applies more generally to the current context of convex, nonincreasing sets $\mathcal{C}$. We do not reproduce it here.

Let then $t$ be a point at which all functions $D_r$ are differentiable. By Rademacher's theorem, this holds for almost every $t \in [0, T]$. Consider first the case where $x_r(t) > 0$. One then has, for all $h > 0$:

$$\int_t^{t+h} \lambda_r^{\mathrm{PF}}(z_{k'}^{-1} X^{k'}(z_{k'} u)) \, du \to \int_t^{t+h} \lambda_r^{\mathrm{PF}}(x(u)) \, du,$$

in view of (i) Lipschitz continuity of $u \to x(u)$, which entails positivity of $x_r(u)$ on $[t, t+h]$, the continuity property of $\lambda_r^{\mathrm{PF}}$ given in Lemma A.3, and finally by an application of Lebesgue's dominated convergence theorem. Therefore, appealing once more to Lemma A.3, the derivative of function $D_r$ at $u$ must coincide with $\lambda_r^{\mathrm{PF}}(x(t))$.

Consider now the case where $x_r(t) = 0$. Clearly, by Fatou's lemma, for all $h > 0$, one has

$$\limsup_{k' \to \infty} \int_t^{t+h} \lambda_r^{\mathrm{PF}}(z_{k'}^{-1} X^{k'}(z_{k'} u)) \, du \leq \int_t^{t+h} \limsup_{y \to x(u)} \lambda_r^{\mathrm{PF}}(y) \, du.$$

On the other hand, the function $x \to \limsup_{y \to x} \lambda_r^{\mathrm{PF}}(y)$ is upper semi-continuous, and thus it follows that

$$\limsup_{u \to t} \left[ \limsup_{y \to x(u)} \lambda_r^{\mathrm{PF}}(y) \right] \leq \limsup_{y \to x(t)} \lambda_r^{\mathrm{PF}}(y).$$

This readily implies that, necessarily, the derivative of $u \to D_r(u)$ at $t$ must lie in the interval $[0, \limsup_{y \to x(t)} \lambda_r^{\mathrm{PF}}(y)]$.

We have thus shown that for any interval $[0, T]$, with probability 1, from any subsequence one can extract a further subsequence $k'$ along which the rescaled process $z_{k'}^{-1} X'k(z_{k'} \cdot)$ converges to a fluid trajectory, uniformly on $[0, T]$. That is to say, almost surely, the accumulation points of the rescaled trajectories consist only of fluid trajectories. This is in contradiction with (31), and the result of Theorem 2 then follows.

## APPENDIX B: PROOF OF LEMMA 5

*Proof of part* (i). Let $u \to x(u)$ denote a fluid trajectory. Let $t > 0$ be a point at which all the associated functions $D_r$ are differentiable. For notational convenience, write $x$ for $x(t)$, $I$ for $I(t)$, $\bar{I}$ for $\mathcal{R} \setminus I$, and $\lambda_r$ for $\lambda_r^{\text{PF}}(x(t))$. For $r \in \bar{I}$, Theorem 2 establishes that $\dot{D}_r(t) = \lambda_r$. Let $d_r$ denote the derivative $\dot{D}_r(t)$ for $r \in I$. As the trajectories $u \to x_r(u)$ are constrained to be nonnegative, necessarily one has for all $r \in I$:

$$\dot{x}_r(t) = 0 = \overline{\nu}_r + \sum_{s \notin I} \mu_s p_{sr} \lambda_s + \sum_{s \in I} \mu_s p_{sr} d_s - \mu_r d_r.$$

This can be written in matrix form as

$$(\mu d)_I = \overline{\nu}_I + (P^T)_{I\bar{I}}(\mu \lambda)_{\bar{I}} + (P^T)_{II}(\mu d)_I,$$

where $(x)_J$ denotes the vector with entries $j \in J$, and $(M)_{IJ}$ denotes the matrix with entries $M_{ij}$, $i \in I$, $j \in J$. The matrix $(P^T)_{II}$ has a spectral radius strictly less than 1, for otherwise the original routing matrix $P$ would have a spectral radius of at least 1. Thus there exists a unique solution $(\mu d)_I$ to the above equation, given by

$$(\mu d)_I = (I - (P^T)_{II})^{-1}[\overline{\nu}_I + (P^T)_{I\bar{I}}(\mu \lambda)_{\bar{I}}].$$

In view of this expression, for $r \in \bar{I}$, the time derivatives $\dot{x}_r(t)$ can be written as

$$\begin{aligned}
(\dot{x})_{\bar{I}} &= (\overline{\nu})_{\bar{I}} + (P^T)_{\bar{I}\bar{I}}(\mu \lambda)_{\bar{I}} + (P^T)_{\bar{I}I}(\mu d)_I \\
&= (\overline{\nu})_{\bar{I}} + (P^T)_{\bar{I}I}(I - (P^T)_{II})^{-1}\overline{\nu}_I \\
&\quad + [(P^T)_{\bar{I}\bar{I}} + (P^T)_{\bar{I}I}(I - (P^T)_{II})^{-1}(P^T)_{I\bar{I}}](\mu \lambda)_{\bar{I}} \\
&= \tilde{\nu} + \tilde{P}^T(\mu \lambda)_I,
\end{aligned}$$

where we have introduced the notation

$$\tilde{\nu} = (\overline{\nu})_{\bar{I}} + (P^T)_{\bar{I}I}(I - (P^T)_{II})^{-1}\overline{\nu}_I,$$
$$\tilde{P}^T = (P^T)_{\bar{I}\bar{I}} + (P^T)_{\bar{I}I}(I - (P^T)_{II})^{-1}(P^T)_{I\bar{I}}.$$

The modified routing probability $\tilde{p}_{rs}$ can be interpreted as the probability that, in a Markov chain on $\mathcal{R}$ started at $r \in \bar{I}$, evolving according to the original routing probabilities $p_{ij}$ (which may become absorbed outside the set $\mathcal{R}$), the next visit to the set $\bar{I}$ is precisely to state $s$. That is to say, $\tilde{p}$ capture the transition probability in the original chain, after removing all excursions to the set $I$. This interpretation allows to deduce at once that the modified routing probability matrix $\tilde{P}$ is sub-stochastic and with spectral radius strictly less than 1 whenever $P$ is so.

It remains to establish the identity

$$(\nu)_{\bar{I}} = (I - \tilde{P}^T)^{-1}\tilde{\nu}.$$

Again, this can be established from a probabilistic interpretation. Assume without loss of generality (by joint rescaling) that the vector $\overline{\nu}$'s entries sum to 1. Then $\nu_r$ can be interpreted as the average number of visits to state $r$ in the Markov chain, with transition probabilities $p$, assuming that the initial distribution is specified by $\overline{\nu}$. It is readily verified that $\tilde{\nu}$ then represents the distribution of the first visit to $\bar{I}$ which is also the initial distribution of the chain where excursions to $I$ are removed. The mean number of visits to states $r \in \bar{I}$ is the same with or without removal of excursions into $I$, hence the desired identity holds.

*Proof of part* (ii). Let us establish (17). In view of Lemma 1, one has

$$(34) \quad L(x(t+h)) - L(x(t)) = \sum_{r=1}^{|\mathcal{R}|} \int_{x_r(t)}^{x_r(t+h)} [\gamma_r^{\mathrm{PF}}(y^r(u)) - \log(\rho_r)]\, du,$$

where the vector $y^r(u)$ is defined as

$$y_s^r(u) = \begin{cases} x_s(t+h), & s < r, \\ u, & s = r, \\ x_s(t), & s > r. \end{cases}$$

At a point $t$ where the fluid trajectories are differentiable, one thus has, by the continuity of functions $\gamma_r^{\mathrm{PF}}$ at $x(t)$ for $r \in \bar{I}$, which follows from Lemma A.3,

$$(35) \quad \lim_{h \to 0} \frac{1}{h} \int_{x_r(t)}^{x_r(t+h)} [\gamma_r^{\mathrm{PF}}(y^r(u)) - \log(\rho_u)]\, du = \dot{x}_r(t)[\gamma_r^{\mathrm{PF}}(x(t)) - \log(\rho_r)].$$

For $r \in I$ and $h > 0$, one has the evaluation

$$\int_{x_r(t)}^{x_r(t+h)} \gamma_r^{\mathrm{PF}}(y^r(u))\, du \leq \sup_{y \in \mathbb{R}_+^{\mathcal{R}}} (\gamma_r^{\mathrm{PF}}(y))(x_r(t+h) - x_r(t)).$$

Indeed, this holds because $x_r(t+h) - x_r(t) \geq 0$, which holds in turn because $x_r(t) = 0$ for $r \in I$, and $x_r(t+h) \geq 0$. This inequality entails that

$$\limsup_{h \searrow 0} \frac{1}{h} \int_{x_r(t)}^{x_r(t+h)} \gamma_r^{\mathrm{PF}}(y^r(u))\, du \leq \sup_{y \in \mathbb{R}_+^{\mathcal{R}}} (\gamma_r^{\mathrm{PF}}(y))\dot{x}_r(t) = 0,$$

where boundedness from above of $\gamma_r^{\mathrm{PF}}$ has been used. This inequality, together with (35) and (34) establish (17).

*Proof of part* (iii). We finally prove (18). To this end, we establish upper bounds on each of the terms in the right-hand side of (34). Note that, for all $t > 0$, and $r \in \bar{I}(x(t))$, the functions $D_r$ are differentiable at $t$, with derivative $\lambda_r^{\mathrm{PF}}(x(t))$. Also, for all $s \in \mathcal{R}$, the functions $D_s$ are nondecreasing and Lipschitz with some constant $A$. It thus follows that, for $r \in \bar{I}$:

$$\limsup_{h \searrow 0} \frac{1}{h} \int_{x_r(t)}^{x_r(t+h)} [\gamma_r^{\mathrm{PF}}(y^r(u)) - \log(\rho_r)] \, du$$
$$\leq \sup_{y \in \mathbb{R}_+^{\mathcal{R}}} (\gamma_r^{\mathrm{PF}}(y) - \log(\rho_r))(\overline{\nu}_r + |\mathcal{R}|A)$$
$$- \mu_r \log(\rho_r) \lambda_r^{\mathrm{PF}}(x(t)) - \mu_r \gamma_r^{\mathrm{PF}}(x(t)) \lambda_r^{\mathrm{PF}}(x(t)).$$

The first two terms in the right-hand side are bounded from above, and the last term is uniformly bounded, since $\lambda_r^{\mathrm{PF}}(x(t)) = \exp(\gamma_r^{\mathrm{PF}}(x(t)))$ and the function $u \to ue^u$ is bounded on a range $(-\infty, A]$.

It remains to consider the case where $r \in I$. One then has

$$\int_{x_r(t)}^{x_r(t+h)} [\gamma_r^{\mathrm{PF}}(y^r(u)) - \log(\rho_r)] \, du$$
$$\leq \sup_{y \in \mathbb{R}_+^{\mathcal{R}}} (\gamma_r^{\mathrm{PF}}(y) - \log(\rho_r))[x_r(t+h) - x_r(t)]$$
$$\leq \sup_{y \in \mathbb{R}_+^{\mathcal{R}}} (\gamma_r^{\mathrm{PF}}(y) - \log(\rho_r))hA,$$

where $A$ is a Lipschitz constant for $u \to x_r(u)$, and nonnegativity of the fluid trajectories has been used. These last two upper bounds together combine to give (18).

## APPENDIX C: PROOF OF LEMMA 6

It follows from Proposition 3, page 21 in [8] (see also [18]) that a continuous function $f$ verifying assumption (20) is such that

(36) $$s < t \Rightarrow f(t) - f(s) \leq (t-s)A.$$

Define now the *increasing variation* $V_f^+(t)$ as the supremum over partitions $\tau_0 = 0 < \tau_1 < \cdots < \tau_m = t$ of the sum

$$\sum_{i=0}^{m-1} (f(\tau_{i+1}) - f(\tau_i))^+.$$

In view of (36), it follows that

$$V_f^+(t) \leq At, \qquad t \in [0, T].$$

It is easily shown that for all $s < t \in [0, T]$, one has

$$V_f^+(t) - V_f^+(s) = \sup_{\tau_0 = s < \cdots < \tau_n = t} \sum_{i=0}^{n-1} (f(\tau_{i+1}) - f(\tau_i))^+,$$

where the supremum is taken over all finite partitions $\tau_0 = s < \cdots < \tau_n = t$. This together with (36) implies that

$$0 \leq V_f^+(t) - V_f^+(s) \leq A(t-s), \qquad s < t \in [0, T].$$

Moreover, if one defines $V_f^-(t)$ as

$$V_f^-(t) := V_f^+(t) - f(t),$$

one readily sees that $u \to V_f^-(u)$ is a nondecreasing function. One may associate a nonnegative measure $\mu^-$ on $[0, T]$ to $V_f^-$ by setting

$$\mu^-([0,t]) = V_f^-(t+) - V_f^-(0).$$

By Radon–Nykodim's theorem, this measure can further be decomposed into a measure that is absolutely continuous with respect to Lebesgue measure, whose density we shall denote by $g^-(t)$, and into a measure $\nu^-$ that is supported by a set $F$ of null Lebesgue measure.

By Rademacher's theorem, the Lipschitz-continuous function $V_f^+$ is almost everywhere differentiable; denote its derivative by $g^+(t)$. Thus, the function $f$ is differentiable almost everywhere, with derivative $g^+(t) - g^-(t)$. Moreover, by condition (19), for almost every $t$, it holds that

$$g^+(t) - g^-(t) \leq -\varepsilon.$$

To conclude, for $s < t < T$, write

$$f(t) - f(s) \leq \int_s^t (g^+(u) - g^-(u)) \, du - \nu^-((s, t))$$
$$\leq -\varepsilon(t-s),$$

which is the announced result.

## APPENDIX D: PROOF OF LEMMA 7

PROOF. Using the notation $x^\pm = \max(0, \pm x)$, note that the factor of $p_{rs}^{(n)}$ in (23) reads

$$(e^{u_s} - 1) \left[ u_s - \sum_{\ell \in \mathcal{R}} p_{s\ell} u_\ell \right] = [(e^{u_s} - 1)^+ - (e^{u_s} - 1)^-] \times \cdots$$

$$\times \left[ u_s^+ - u_s^- - \sum_{\ell \in \mathcal{R}} p_{s\ell} u_\ell^+ + \sum_{\ell \in \mathcal{R}} p_{s\ell} u_\ell^- \right]$$

$$= (e^{u_s} - 1)^+ \left[ u_s^+ - \sum_{\ell \in \mathcal{R}} p_{s\ell} u_\ell^+ \right]$$

$$+ (e^{u_s} - 1)^- \left[ u_s^- - \sum_{\ell \in \mathcal{R}} p_{s\ell} u_\ell^- \right]$$

$$+ (e^{u_s} - 1)^+ \sum_{\ell \in \mathcal{R}} p_{s\ell} u_\ell^-$$

$$+ (e^{u_s} - 1)^- \sum_{\ell \in \mathcal{R}} p_{s\ell} u_\ell^+.$$

In order to obtain the above expansion, we have used the fact that $(e^{u_s}-1)^{\pm} u_s^{\mp} = 0$. Note that the last two terms in this expansion are nonnegative. Note also that the first two terms both read

$$(e^{v_s} - 1)\left[ v_s - \sum_{\ell \in \mathcal{R}} p_{s\ell} v_\ell \right]$$

for adequate choices of $v_s$, namely $v_s = u_s^+$ for the first term, and $v_s = -u_s^-$ for the second term. This establishes claim (ii) of the lemma. This further implies that, in order to prove claims (i) and (iii) of the lemma, it is sufficient to restrict attention to the case where the $u_s$ all have the same sign, which we now assume.

Introduce the notation

$$\begin{cases} M(s,\ell) := \sum_{n \geq 0} p_{rs}^{(n)} p_{s\ell}, \\ N(s,\ell) := M(s,\ell) + \mathbf{1}_{\ell=r} \sum_{n \geq 0} p_{rs}^{(n)} \left( 1 - \sum_{k \in \mathcal{R}} p_{sk} \right) \\ \phantom{N(s,\ell)} = \sum_{n \geq 0} p_{rs}^{(n)} \left( p_{s\ell} + \mathbf{1}_{\ell=r}\left( 1 - \sum_{k \in \mathcal{R}} p_{sk} \right) \right). \end{cases}$$

Condition (23) thus reads

$$\sum_{s,\ell \in \mathcal{R}} M(s,\ell)(e^{u_s} - 1)u_\ell \leq \sum_{s,\ell \in \mathcal{R}} N(s,\ell)(e^{u_s} - 1)u_s.$$

Let us now show that this last condition is satisfied. Note first that it is enough to prove the same inequality, with $N$ instead of $M$ in the left-hand side, since

$$\sum_{s,\ell \in \mathcal{R}} N(s,\ell)(e^{u_s} - 1)u_\ell - \sum_{s,\ell \in \mathcal{R}} M(s,\ell)(e^{u_s} - 1)u_\ell$$

(37)
$$= \sum_{s \in \mathcal{R}} \sum_{n \geq 0} p_{rs}^{(n)} \left(1 - \sum_{k \in \mathcal{R}} p_{sk}\right)(e^{u_s} - 1)u_r,$$

and this difference is indeed nonnegative under the current assumption that the $u_s$ all have the same sign.

We thus need to show that

(38) $$\sum_{s,\ell \in \mathcal{R}} N(s,\ell)(e^{u_s} - 1)u_\ell \leq \sum_{s,\ell \in \mathcal{R}} N(s,\ell)(e^{u_s} - 1)u_s.$$

Note now that the marginals of the measure $N(\cdot,\cdot)$ coincide. Indeed,

$$\sum_{\ell \in \mathcal{R}} N(s,\ell) = \sum_{n \geq 0} p_{rs}^{(n)},$$

$$\sum_{s \in \mathcal{R}} N(s,\ell) = \sum_{n \geq 0} p_{r\ell}^{(n+1)} + \mathbf{1}_{\ell = r} \sum_{n \geq 0} \sum_{s \in \mathcal{R}} (p_{rs}^{(n)} - p_{rs}^{(n+1)})$$

$$= \sum_{n \geq 0} p_{r\ell}^{(n)}.$$

Thus, after renormalization of both sides of (38) by the total mass of the measure $N$, it equivalently reads

(39) $$\mathbf{E}[(e^U - 1)V] \leq \mathbf{E}[(e^U - 1)U],$$

where the random variables $U$, $V$ have the same distributions. An inequality due to Hoeffding [12] (see also [28] and [7] for more easily accessible references) states that, given two random variables $U$, $V$ with identical distributions, for any two nondecreasing functions $f, g : \mathbb{R} \to \mathbb{R}$ such that $f(U)$ and $g(U)$ have finite variances, one has

$$\mathbf{E}[f(U)g(V)] \leq \mathbf{E}[f(U)g(U)].$$

Note that the inequality we need to prove is of that form, with as nondecreasing functions $f(U) = U$ and $g(U) = e^U - 1$. Finiteness of variances is trivially satisfied as the random variables $U$ take only finitely many values. This concludes the proof of the first claim of the lemma.

Let us now show that, in order to have equality in (23), all $u_s$ such that $\sum_{n \geq 0} p_{rs}^{(n)} > 0$ must be zero. Equality in (23) implies equality in (39). However, the latter holds if and only the distributions of $(f(U), g(V))$ and $(f(U), g(U))$ coincide; see, for example, [28]. As the functions $f$, $g$ are strictly increasing, this in turn holds if the distributions of $(U, V)$ and $(U, U)$ coincide. This means that we can partition the set $\mathcal{R}$ such that on each subset of the partition, the $u_s$ are constant, and for $s, \ell$ in two different subsets of the partition, $N(s,\ell) = 0$. Thus, for all $s$ such that $\sum_{n \geq 0} p_{rs}^{(n)} > 0$, one must

32    L. MASSOULIÉ

have $u_r = u_s$. This is needed to ensure equality in (38) . However, in order to ensure equality in (23), the right-hand side of (37) must also be zero, which, using the fact that all $u_s$ coincide, also reads

$$r \in \bar{I} \Rightarrow w_r = \gamma_r^{\mathrm{PF}}(x) = \log(\lambda_r^{\mathrm{PF}}(x))u_r(e^{u_r} - 1) \sum_{s \in \mathcal{R}} p_{rs}^{(n)} \left(1 - \sum_{k \in \mathcal{R}} p_{sk}\right) = 0.$$

Equivalently, one must have $u_r(\exp(u_r) - 1) = 0$, that is, $u_r = 0$, which concludes the proof of the lemma. $\square$

REMARK A.1.   Note that the statement of Lemma 7 remains true if we replace the terms $[\exp(u_s) - 1]$ by $f(u_s)$ in (23), where $f$ is any strictly increasing function such that $f(0) = 0$.

Acknowledgments.**Acknowledgments.**  The author is very grateful to the editor and reviewers, whose comments greatly helped improve this work. In particular, they pointed to a serious gap in an earlier version. Theorem 2, Lemmas 5 and 6 have been established and included in the present version so as to fill the gap in the previous treatment.

The author would also like to acknowledge stimulating discussions related to this work with Maury Bramson, A. J. Ganesh, Frank Kelly and Ruth Williams.

## REFERENCES

[1] ASMUSSEN, S. (2003). *Applied Probability and Queues.* Springer, New York. MR1978607
[2] BERTSEKAS, D. and GALLAGER, R. (1992). *Data Networks.* Prentice Hall, New York.
[3] BONALD, T. and MASSOULIÉ, L. (2001). Impact of fairness on Internet performance. In *Proceedings of ACM Sigmetrics* 82–91.
[4] BONALD, T., MASSOULIÉ, L., PROUTIÈRE, A. and VIRTAMO, J. (2006). A queueing analysis of max-min fairness, proportional fairness and balanced fairness. *Queueing System Theory Appl.* **53** 65–84. MR2230014
[5] BONALD, T. and PROUTIÈRE, A. (2003). Insensitive bandwidth sharing in data networks. *Queueing Systems* **44** 69–100. MR1989867
[6] BRAMSON, M. (2005). Communication at the INFORMS Applied Probability Conference, Ottawa.
[7] BRATLEY, P., FOX, B. L. and SCHRAGE, L. E. (1986). *A Guide to Simulation.* Springer, New York.
[8] CERAGIOLI, F. M. (1999). Discontinuous ordinary differential equations and stabilization. Ph.D. dissertation, Univ. degli Studi di Firenze.
[9] DAI, J. G. (1995). On positive Harris recurrence of multiclass queuing networks: A unified approach via fluid limit models. *Ann. Appl. Probab.* **5** 49–77. MR1325041
[10] DE VECIANA, G., LEE, T. J. and KONSTANTOPOULOS, T. (1999). Stability and performance analysis of networks supporting services with rate control—could the Internet be unstable? In *Proc. IEEE Infocom'99.*

[11] FALCONER, K. (1990). *Fractal Geometry. Mathematical Foundations and Applications*. Wiley, Chichester. MR1102677
[12] HOEFFDING, W. (1940). Masstabinvariante korrelationstheorie. *Schriften des Mathematischen Instituts und des Instituts für Angewandte Mathematik der Universitat Berlin* **5** 179–233.
[13] KELLY, F. (1997). Charging and rate control for elastic traffic. *Eur. Trans. Telecommun.* **8** 33–37.
[14] KELLY, F., MAULLOO, A. and TAN, D. (1998). Rate control in communication networks: Shadow prices, proportional fairness and stability. *J. Oper. Res. Soc.* **49** 237–252.
[15] KELLY, F. P. and WILLIAMS, R. J. (2004). Fluid model for a network operating under a fair bandwidth-sharing policy. *Ann. Appl. Probab.* **14** 1055–1083. MR2071416
[16] KEY, P. and MASSOULIÉ, L. (2006). Fluid models of integrated traffic and multipath routing. *Queueing Systems Theory Appl.* **53** 85–98. MR2230015
[17] MAZUMDAR, R., MASON, L. and DOULIGERIS, C. (1991). Fairness in network optimal flow control: Optimality of product forms. *IEEE Trans. Commun.* **39** 775–782.
[18] MCSHANE, E. J. (1947). *Integration*. Princeton Univ. Press. MR0082536
[19] MO, J. and WALRAND, J. (2000). Fair end-to-end window-based congestion control. *IEEE/ACM Transactions on Networking* **8** 556–567.
[20] NASH, J. (1950). The bargaining problem. *Econometrica* **18** 155–162. MR0035977
[21] ROBERT, P. (2003). *Stochastic Networks and Queues*. Springer, Berlin. MR1996883
[22] ROCKAFELLAR, R. T. (1970). *Convex Analysis*. Princeton Univ. Press. MR0274683
[23] ROCKAFELLAR, R. T. and WETS, R. J. B. (2004). *Variational Analysis*. Springer, Berlin. MR1491362
[24] RYBKO, A. N. and STOLYAR, A. L. (1992). Ergodicity of stochastic processes describing the operations of open queueing networks. *Problemy Pederachi Informatsii* **28** 2–26. MR1189331
[25] SCHASSBERGER, R. (1986). Two remarks on insensitive stochastic models. *Adv. in Appl. Probab.* **18** 791–814. MR0857330
[26] SHORAK, G. and WELLNER, J. (1986). *Empirical Processes with Applications to Statistics*. Wiley, New York. MR0838963
[27] STEFANESCU, A. and STEFANESCU, M. W. (1984). The arbitrated solution for multiobjective convex programming. *Rev. Roum. Math. Pure Appl.* **29** 593–598. MR0759520
[28] WHITT, W. (1976). Bivariate distributions with given marginals. *Ann. Statist.* **4** 1280–1289. MR0426099
[29] YE, H. Q. (2003). Stability of data networks under an optimization-based bandwidth allocation. *IEEE Trans. Automat. Control* **48** 1238–1242. MR1988097
[30] YE, H. Q., OU, J. H. and YUAN, X. M. (2005). Stability of data networks: Stationary and bursty models. *Oper. Res.* **53** 107–125. MR2131101

PARIS RESEARCH LAB
THOMSON
46 QUAI A. LE GALLO
92648 BOULOGNE CEDEX
FRANCE
E-MAIL: laurent.massoulie@thomson.net